\documentclass[11pt,a4paper]{article}
\usepackage{amsmath,amssymb}
\newcommand{\n}{\nabla}

\newcommand{\p}{\partial}

\newtheorem{definition}{Definition}
\newtheorem{theorem}{Theorem}
\newtheorem{proposition}{Proposition}

\newtheorem{remarka}{Remark}

\newtheorem{lemme}{Lemma}

\title{Remarks on global controllability for the shallow-water system  with two control forces}
\author{Abdelmalek Drici \footnote{UPMC Univ  Paris 06, UMR 7598 Laboratoire Jacques-Louis Lions, Paris, F-75005 France ;
  CNRS, UMR 7598 LJLL, Paris, F-75005 France}, Boris Haspot \thanks{Ceremade UMR CNRS 7534
Universit\'e de Paris  Dauphine,
Place du MarŽchal DeLattre De Tassigny
75775 PARIS CEDEX 16 , haspot@ceremade.dauphine.fr }}

\date{}
\begin{document}
%\tableofcontents

\maketitle
%VOIR COMMENT UTILISER LA METHODE A CORON POUR CREER ARTIFICIELLEMENT DU DAMPING!!!
\begin{abstract}
In this paper we deal with the compressible Navier-Stokes equations with a friction term in one dimension on an interval. We study the exact controllability properties of this equation with general initial condition when the boundary control is acting at both endpoints of the interval. Inspired by the work of Guerrero and Imanuvilov in \cite{GI} on the viscous Burger equation, we prove by choosing irrotational data and using the notion of effective velocity developed in \cite{cpde,cras} that the exact global controllability result does not hold for any time $T$. 
\end{abstract}
\section{Introduction}
We consider the viscous shallow water model with friction term. This model is also called by the french community the Saint-Venant equations and is generally used in oceanography. Indeed it allows to model vertically averaged flows  in terms of the horizontal mean velocity field $u$ and the depth variation $h$. In the rotating framework, the model is described by the following system:
\begin{equation}
\begin{cases}
\begin{aligned}
&\frac{\p}{\p t}h+\p_{x}(h u)=0,\;\;(t,x)\in Q=(0,T)\times (0,1)\\
&\frac{\p}{\p t}(h u)+\p_{x}(h
|u|^{2})-\p_{x}(\mu h\p_{x}u)+\frac{\p_{x} h}{Fr^{2}}+r h u=0,\;\;(t,x)\in Q=(0,T)\times (0,1)\\
&u(t,0)=v_{1}(t),\;\;u(t,1)=v_{2}(t),\;\;t\in(0,T),
%&u(0,x)=u_{0}(x),\;h(0,x)=h_{0}(x).
\end{aligned}
\end{cases}
\label{3systeme}
\end{equation}
$Fr > 0$ denotes the Froude number and we set $r=\frac{1}{\mu Fr^{2}}$. System (\ref{3systeme}) is supplemented with initial conditions
\begin{equation}
h_{/t=0}= h_{0}, u_{/t=0}=u_{0}.
\label{data}
\end{equation}
%This model is derived from the three-dimensional Navier-Stokes equations with free surface, where the normal stress is determined from the air pressure and capillary effects. 
This model is derived from the compressible Navier-Stokes equations and the turbulent regime ($r\geq 0$) is obtained from the friction condition on the bottom, see \cite{P}. $\mu$ is the viscosity coefficient and verifies $\mu>0$ .\\
Here, $T > 0$ is a given final time and $v_{1}(t)$ and $v_{2}(t)$ are control functions which are acting over our system at both endpoints of the segment $(0,1)$. We are now going to assume that
the controls $v_{1}$ and $v_{2}$ belong to the space $H^{3/4}(0,T)$ and the initial data  $(h_{0},u_{0})$  are in $H^{2}(0,1)\times H^{1}(0,1)$ with $h_{0}\geq c>0$. Furthermore we suppose the following compatibility assumption:
%VOIR LA PARTIE EXISTENCE, ON PEUT LE FAIRE A L'INVERSE PAR LES EQUIVALENCES DE GUERRERO.\\
%Furthermore, $(h_{0},u_{0})$ are the initial condition which is supposed to be in $H^{1}(0,1)$. In the sequel, we will suppose that our control functions $v_{1}$ and $v_{2}$ belong to the space $H^{3/4}(0,T)$ and they satisfy the compatibility conditions:
\begin{equation}
v_{1}(0) = u_{0}(0)\;\;\mbox{and}\;\; v_{2}(0) = u_{0}(1). 
\label{2}
\end{equation}
Under these assumptions, it is classical to see that there exists a solution $(h,u)$ of system (\ref{3systeme}) which belongs for a $T>0$ to the space
$X_{T}=L^{\infty}(0,T;H^{2}(0,1))\times\big(L^{2}(0,T;H^{2}(0,1))\cap H^{1}(0,T;L^{2}(0,1))\big)$ and  a constant $C>0$ such that:
\begin{equation}
\|(h,u)\|_{X_{T}} \leq C (\|u_{0}\|_{H^{1}(0,1)} +\|h_{0}\|_{H^{2}(0,1)} +\|v_{1}\|_{H^{3/4}(0,T)} +\|v_{2}\|_{H^{3/4}(0,T)} ),
\label{3}
\end{equation}
(see, for instance, \cite{11}). We would like to mention that we can obtain global solution provided that we choose $u_{0}=-\mu\p_{x}\ln h_{0}$ (we shall detail more precisely in the sequel this point which will play a crucial role in your results, we refer also to \cite{cras}).\\
In the sequel we are interesting in proving some negative results about the global exact controllability of the system (\ref{3systeme}). Before giving some elements of the proof we would like to recall the notion of global exact controllability property for system (\ref{3systeme}).
\begin{definition}
\label{control}
The system  (\ref{3systeme}) will be said global exact controllable if for any $(h_{0},u_{0})\in (H^{2}(0, 1)\times H^{1}(0, 1))$ with $h_{0}\geq c>0$, and for any $(h_{1},u_{1})\in (H^{2}(0, 1)\times H^{1}(0, 1))$ with $h_{1}\geq c_{1}>0$,  for any time $T>0$ there exist controls $v_{1}\in H^{ 3/4} (0, T )$ and $v_{2}\in H^{ 3/4}  (0, T )$ such that the corresponding solution of (\ref{3systeme})  satisfies $u(T , x) = u_{1}(x)$ in $(0, 1)$ and $h(T,x)= h_{1}(x)$ in $(0, 1)$.%, we will refer to this problem as \textbf{the exact null controllability}.
\end{definition}
\begin{remarka}
In our previous definition, we impose only two controls on the system (\ref{3systeme}) which is natural inasmuch as the result of existence of strong solution for the system (\ref{3systeme}) do not admit boundary condition on the height $h$.% ICI ON NE MET QUE DEUX CONTROLES ALORS QU'ON POURRAIT EN METTRE PLUS (QUATRE EN TOUT)!!!!
\end{remarka}
Let us recall what is known about the controllability of the compressible Navier-Stokes equations, as far as we know there exists very few results. Indeed one of the main difficulty consists in the non controllability of the linear system associated to (\ref{3systeme}) around the stable equilibrium $(1,0)$. This has been observed indirectly in \cite{RR} and more recently in \cite{DH}. By developing a moment method in one dimension, we can observe that the operator associated to the linear system of (\ref{3systeme}) admits accumulation point on the eigenvalues, this is due in particular to the behavior of the height which generates a damping effect. However in \cite{DH}, we check that the linearized system is approximate controllable. To overcome this difficulty, it would be natural to use the return method developed by Coron (see \cite{Cor1,4,5}), however the linearized system around the trajectory of \cite{Cor1} remains non controllable. It is one of the main reason why the problem of exact controllability for the compressible Navier-Stokes equations is so difficult. However very recently in a very interesting work Ervedoza, Glass, Guerrero and  Puel in \cite{EGGP} have obtained the local exact controllability for the 1-D compressible Navier-Stokes equation by using Carleman estimates taking into account the flow associated to the transport equation on the height.\\
The study of the exact local incompressible Navier-Stokes equations has been investigated by many authors, we would like in particular recall the results of Fern\'andez-Cara et al in \cite{FG1,FG2}. Global exact controllability in two dimension has been obtained by Coron and  Fursikov for a manifold without boundary by combining some result of global approximate control obtained via the return method and some results of local exact controllability via Carleman estimates.
\\
The goal of this paper is to show that the compressible shallow water system with friction term in not global exact controllable and not even  global approximate controllable.  To do this, we are going constructed special solution of the system (\ref{3systeme}) which are not controllable (they are also called \textit{quasi solution} in \cite{ras} when there is no friction term). Here
compared with compressible Navier-Stokes equation we also take into account the friction term $r h u$ physically justified to model the friction condition on the bottom of the ocean ( see \cite{P}). Naturally on a mathematical point of view, the term $r h u$ does not add any difficulties for obtaining global weak solution or global strong solution with small initial data, however in our study this term shall turn out to be essential in order to cancel out the coupling between the height $h$ and the velocity $u$. Initially in \cite{cpde} we  obtain the existence of global strong solution for Korteweg system with friction term with large initial data (at least for the irrotational part of the initial data) in dimension $N\geq 2$. Indeed coupled with the pressure term, the friction term shall introduce a damping effect on a new unknown called \textit{effective velocity} and introduced in \cite{Hprepa1}.
 Roughly speaking this friction term allows to cancel out the coupling between the height $h$ and the velocity $u$ (we can check in particular that $(h,-\mu\n\ln h)$ is a particular solution of the system where $h$ verifies a heat equation), it is one of the main difficulty in order to obtain the existence of global strong solution for compressible Navier-Stokes equation (see \cite{arma}, in this case the coupling is between the velocity and the density). Indeed in this last case, it is difficult to obtain a damping effect on the density and the coupling between the velocity and the pressure terms impose a smallness condition on the initial density.\\
 This last result was extended for the shallow water system without friction in  \cite{cras}, in this paper we introduce the notion of quasi-solution $(h_{1},-\mu\n\ln h_{1})$ ( where $h_{1}$ verifies a heat equation) which is a particular solution of the shallow-water system  when the pressure is considered null. We consider then the pressure term as a perturbation by taking profit of the regularizing effect on $h_{1}$ which induces smallness
condition on the pressure in high frequencies.\\
In this paper we are going to verify that $(h,-\mu\p_{x}\ln h)$ is a particular solution of (\ref{3systeme}). Via Cole-Hopf formula we verifies that the global controllability problem of (\ref{3systeme}) comes to deal in reality with the viscous Burger equation studied in particular by Guerrero and Imanuvilov in \cite{GI}. They  prove that this equation is not global exact controllable by using some comparison argument on the heat equation. We shall follow in our proof the argument of \cite{GI} which shall be an angular stone in our result.\\
Before giving our main result, let us briefly mention some previous works which deal with the exact controllability for the viscous and non viscous Burgers equation in order to understand better this equation. 
\subsection*{Exact controllability for the non viscous Burger equations}
%We start with recalling some results about the non-viscous Burgers equation. As far as we know, only two works have been dedicated to this issue.\\
In \cite{10},   Horsin shows that by using one boundary control, the Burgers equation can be driven from the null initial condition to a constant final state $M$ in a time $T\geq 1/M$. To do this, he  uses is the so-called return method introduced by  Coron (see \cite{5} for the case of the global approximate control of the Navier-Stokes equation with slip boundary condition). More general results are obtained in \cite{1} when we consider general scalar non-linear conservation laws.% with $C_{2}$ strictly convex flux functions when starting from a null initial data. In particular, they deduce that as long as the controllability time is lower than $1/M$ ($M > 0$ constant), we do not have exact controllability to the state $M$.\\
\subsection*{Exact controllability for the  viscous Burger equations}
\subsubsection*{With one boundary control}
First, in \cite{6}, the author proves that  the approximate controllability to some target states does not arise for the solutions of the viscous Burgers equation with one boundary control. Indeed he proves the existence of some invariant, more precisely that there exists a constant $C > 0$ such that, for every $T > 0$ one has:
$$u(t,x)\leq \frac{C_{0}}{1-x}\;\;\forall (t,x)\in(0,T)\times(0,1). $$
%From this a priori estimate, the author deduces that the approximate controllability to some target states does not hold. We recall now a result from \cite{8}, where the authors prove that we cannot reach an arbitrary target function in arbitrary time with the help of one control force. Precisely, they deduced the following estimate: for each $N > 5$, there
%exists a constant $C_{1}(N) > 0$ such that:
%$$\frac{d}{dt}\int^{b}_{0}(b-x)^{N}u_{+}^{4}(t,x)dx<C_{1}b^{N-5}$$
%where $b$ is the lower endpoint where the control function is supported (for instance, $b = 1/2$ if the control function is acting in $(1/2,1)$). Here, $u_{+}(t,x) = \max{u(t,x),0}$ is the positive part of $u$. From this a priori estimate, one can deduce that we can not get close to some open set of target functions in $L^(0,1)$.\\
Next, in \cite{7}, the authors prove that the global null controllability does not hold with one control force by using in a crucial way the comparison principle.% Precisely, for any initial condition $u_{0}\in L^{2}(0, 1)$ with $u_{0}\in L^{2}(0,1) = r > 0$, it is proved that there exists a time $T (r) > 0$ such that for any control function, the corresponding solution satisfies
%$$|u(t,\cdot)|\geq C_{2}>0 \;\;\forall t\in(0,T(r))\;\;\mbox{in any open intervalI } I\subset(0,1),$$
%for some $C_{2}> 0$. The proof use crucially the comparison principle.
% Furthermore, this time $T¨$ is proved to be sharp in the sense that there exists a constant $C_{3} > 0$ independent of $r$ such that, if $T > C_{3}T¨ $, then there exists a control function such that the corresponding solution satisfies $u(T , x) = 0$ in $(0, 1)$. The main tool used in this work is the comparison principle.\\
%\\
%For the Burgers equation with two boundary controls it was shown in \cite{8}, that any steady state solution is reachable for a sufficiently large time.
%Finally, in the recent paper \cite{4}, the author proved that with the help of two control forces, we can drive the solution of the Burgers equation with null initial condition to large constant states. More precisely, for any time $T > 0$, it is shown the existence of a constant $C_{4} > 0$ such that for any $C\in\R$ satisfying $|C|\geq C4$ (INTERESSANT DANS MON CAS!!!), there exist two controls $v_{1}(t)$ and $v_{2} (t )$ such that the associated solution to (\ref{3systeme}) with $u_{0}= 0$ satisfies $y (T , á) = C$ in $(0, 1)$. The idea of the proof of this result is based on the HopfÐCole transformation, which leads to a controllability problem for the heat equation.
\subsubsection*{With two boundary controls}
In \cite{GI}, Guerrero and Imanuvilov prove that the viscous Burger equation is not exact global controllable even with two boundary controls. To do this they combine the Cole-Hopf transformation and some comparison principle for the heat equation.\\
In \cite{GG}, Glass and Guerrero show the uniform exact controllability for the viscous Burger equation when the viscosity is small enough. In \cite{IP}, Imanuvilov and Puel obtain the global exact controllability for a type of Burger equations.
To finish we would like to mention the works of  Chapouly in \cite{Cha} who proves the global exact controllability of the viscous Burger equation when we consider three controls: two are the boundary values, one is the right member of the equation and is constant with respect to the space variable. This result is to relate with the paper of  Coron \cite{5} inasmuch as the right control replace in some sense the pressure term of the Navier-Stokes equations. In some sense in \cite{Cha}, Chapouly via the source control forces the system to be incompressible as for Navier-Stokes.\\
In the present paper, we have two main objectives. One concerning the global exact null controllability for small time and the other one concerning the global exact controllability for any time $T > 0$. Both results are of negative nature.\\
As long as the first one is concerned, we prove that there exists a final time $T$ and an initial condition $u_{0}$ such that the solution of (\ref{3systeme}) is far away from zero. That is to say, the global null controllability for the compressible Navier-Stokes equation with friction and with two control forces does not hold. The precise result is given in the following theorem:
\begin{theorem}
 There exists $T > 0$ and $u_{0}=-\mu\p_{x}\ln h_{0}\in H^{1}(0, 1)$ with $h_{0}\geq c>0$ such that, for any control functions $v_{1}\in H^{ 3/4}(0, T )$ and $v_{2}\in H^{3/4} (0, T )$ satisfying the compatibility conditions (\ref{2}), the associated solution $(h,u)\in X_{T}$ to (\ref{3systeme}) satisfies:
\begin{equation}
\|u(T,\cdot)\|_{H^{1}(0,1)}\geq C >0, 
\label{4}
\end{equation}
for some positive constant depending on $T$ and $h_{0}$.
\label{theo1}
\end{theorem}
The second main result is a negative exact controllability result:
\begin{theorem}
 For any $T > 0$, there exists an initial condition $u_{0}-\mu\p_{x}\ln h_{0} \in H^{1} (0, 1)$ and a target function $u_{1}=-\mu\p_{x}\ln h_{1}\in H^{1} (0, 1)$ such that, for any $v_{1}\in H^{3/4} (0, T )$  and $v_{2}\in H^{3/4} (0, T )$ satisfying (\ref{2}), the associated solution $(h,u)\in X_{T}$ to (\ref{3systeme}) satisfies:
 \begin{equation}
 \|u(T,á)-u_{1}(á)\|_{H^{1}(0,1)} \geq C >0, 
  \label{5}
  \end{equation}
for some positive constant $C$ depending on $T$, $h_{0}$ and $h_{1}$.
\label{theo2}
\end{theorem}
\begin{remarka}
Let us observe that we have global existence of solution for the system (\ref{3systeme}) via the fact that $(h,-\mu\p_{x}\ln h)$ is a particular solution.
\end{remarka}
\begin{remarka}
Let us remark that in \cite{5}, J-M Coron proves the global approximate controllability of the Navier-Stokes equations with slip boundary condition by using the return method (in \cite{CF}, the authors obtain the global exact controllability). In \cite{GI}, Guerrero and Imanuvilov prove that the viscous Burger equation in one dimension is not global exactly controllable. We show that it is also the case for the compressible Navier-Stokes equations with friction. This is due in particular to the fact that the flow is compressible, indeed in our case we are working with a specific solution of the form $u=-\mu\p_{x}\ln h$ which is naturally not incompressible but irrotational. In particular it gives an idea why in the case of slip boundary condition we can hope global controllability. Indeed by imposing a Neuman condition on the height which verifies an heat equation in the context of the \textit{quasi solution}, we verify that $u=-\mu\p_{x}\ln h$ is such that $u\cdot n=0$.
\end{remarka}
\begin{remarka}
We recall that all global controllability properties obtained up to now for the incompressible NavierÐStokes system come essentially from the same property for the Euler equation (see the works of Coron \cite{4,9}). In the case of the compressible Navier-Stokes equation we use the properties of the viscous Burger equations and not these of the non viscous Burger equations as for incompressible Navier-Stokes equations. On the other side the structure of \textit{quasi solution} is purely non linear.
\end{remarka}
\begin{remarka}
In a forthcoming paper we are going to prove the same type of results for the shallow water system without friction by using the same type of argument than in \cite{cras} which corresponds to neglect in high frequencies the pressure term.
\end{remarka}
In order to prove these results, in section \ref{section2} we first show the equivalence of the controllability problem for the  system (\ref{3systeme}) and some controllability problem for a one-dimensional linear heat equation with positive boundary controls. To do this, we combine a Cole-Hopf transformation and the introduction of particular solution introduced in \cite{Hprepa1}.
Finally, in the section \ref{section3} we prove both theorems \ref{theo1} and \ref{theo2} by using a method developed by Guerrero and Imanuvilov in \cite{GI}.
 \section{Simplification of the problem}
 \label{section2}
In this section, we show that with particular initial data (it means irrotational initial data on the velocity which depends on $h_{0}$) the exact controllability properties for the system (\ref{3systeme}) (see definition \ref{control}) can be reduced to some controllability assumptions for the heat equation. In the sequel $T > 0$ is a fixed final time.
Let us rewrite the system (\ref{3systeme}):
\begin{equation}
\begin{cases}
\begin{aligned}
&\frac{\p}{\p t}h+\p_{x}(h u)=0,\;\;(t,x)\in Q=(0,T)\times (0,1)\\
&\frac{\p}{\p t}u+u\p_{x}u-\mu\p_{xx}u+\frac{\p_{x} \ln h}{Fr^{2}}+r u=0,\;\;(t,x)\in Q=(0,T)\times (0,1)\\
&u(t,0)=v_{1}(t),\;\;u(t,1)=v_{2}(t),\;\;t\in(0,T),\\
&u(0,x)=u_{0}(x),\;h(0,x)=h_{0}(x).
\end{aligned}
\end{cases}
\label{3systeme1}
\end{equation}
Next we remark that the controllability problem of definition \ref{control} implies the controllability problem for a semilinear parabolic system  (\ref{3systeme2}) with time-dependent controls.%, one acting in the right-hand side of our equation and the other one acting at $x = 1$. 
We are writing the system (\ref{3systeme2})  as follows:
\begin{equation}
\begin{cases}
\begin{aligned}
&\frac{\p}{\p t}h-\mu\p_{xx}h=0,\;\;(t,x)\in Q=(0,T)\times (0,1)\\
&\frac{\p}{\p t}w-\mu\p_{xx}w-\mu|\p_{x}w|^{2}=v^{'}_{5}(t)+v_{3}(t),\;\;(t,x)\in Q=(0,T)\times (0,1)\\%v^{'}_{5}(t)+\\
&w(t,0)=v_{5}(t)
,\;\;w(t,1)=v_{4}(t),\;\;t\in(0,T),\\
&h(t,0)=e^{v_{5}(t)},\;\;h(t,1)=e^{v_{4}(t)},\;\;t\in(0,T),\\
&w(0,x)=\ln h_{0}(x),\;h(0,x)=h_{0}(x).\\
\end{aligned}
\end{cases}
\label{3systeme2}
\end{equation}
\begin{lemme}
The system (\ref{3systeme}) is global exact controllable implies that there exists a solution to the following controllability problem:\\
For any $w_{0},w_{1}\in H^{2}(0,1)$ and $h_{0}=e^{w_{0}}$, $h_{1}=e^{w_{1}}$ with $w_{0}(0)=w_{1}(0)=0$, there exists control $v_{3}\in L^{2}(0,T)$ and $v_{4},v_{5}\in H^{1}(0,T)$ and a solution of (\ref{3systeme2}) such that $w(T,x)=w_{1}(x)$ in $(0,1)$ and $h(T,x)=h_{1}(x)$ with $w\in Y_{T}=L^{2}(0,T;H^{3}(0,1))\cap H^{1}(0,T;H^{1}(0,1))$ and $h=e^{w}$.%\\
%METTRE $v_{5}$
\end{lemme}
{\bf Proof:} It suffices to set:
$$\ln h(t,x)=w(t,x)=-\frac{1}{\mu}\int^{x}_{0}u(t,y)dy+v_{5}(t),\;\;\forall (t,x)\in Q,$$%+v_{5}(t),\;\;\forall (t,x)\in Q,$$
and with:
$$
\begin{aligned}
&v_{3}(t)=-\frac{v_{1}^{2}(t)}{\mu}+u_{x}(t,0),\;v_{4}(t)=-\frac{1}{\mu}\int^{1}_{0}u(t,y)dy+v_{5}(t),\;w_{0}(x)=-\frac{1}{\mu}\int^{x}_{0}u_{0}(y)dy,\\
&w_{0}(x)=-\frac{1}{\mu}\int^{x}_{0}u_{0}(y)dy+v_{5}(0),\;w_{1}(x)=-\frac{1}{\mu}\int^{x}_{0}u_{1}(y)dy+v_{5}(1).
\end{aligned}
$$
From (\ref{3}) and $v_{1}\in H^{\frac{3}{4}}(0,T)$, we obtain that $v_{3}\in L^{2}(0,T)$ and $v_{4},v_{5}\in H^{1}(0,T)$.
%\\
%VOIR LES THEOREMES DE TRACE!!!\\
%\\
%IL S'AGIT DE FAIRE TRES ATTENTION AUX CONDITIONS AUX BORDS SUR RHO ET LA COMPATIBILITE ENTRE LES DEUX EQUATIONS QUAND ON VA REVENIR SUR LA PREMIERE. IL FAUT EN PLUS UTILISER L'UNICITE POUR VERIFIER QUE $w=\ln\rho$. ON REMARQUE ICI QUE PAR INJECTION DE SOBOLEV $w_{1}$ EST DANS $L^{\infty}$ DONC TOUT LES $h_{0}$ SONT ATTEINTS
$\blacksquare$
\\
Let us now prove that the previous controllability result is equivalent to a controllability problem for the system  (\ref{3systeme3}) with two time-dependent controls, one of bilinear nature (multiplying the state function) and the other one acting at $x = 1$. We are writing the system (\ref{3systeme3}):
\begin{equation}
\begin{cases}
\begin{aligned}
&\frac{\p}{\p t}h-\mu\p_{xx}h=0,\;\;(t,x)\in Q=(0,T)\times (0,1)\\
&\frac{\p}{\p t}h-\mu\p_{xx}h-(v^{'}_{5}(t)+v_{3}(t))h=0,\;\;(t,x)\in Q=(0,T)\times (0,1)\\
&h(t,0)=e^{v_{5}(t)},\;\;h(t,1)=e^{v_{4}(t)},\;\;t\in(0,T),\\
&h(0,x)=h_{0}(x).\\
\end{aligned}
\end{cases}
\label{3systeme3}
\end{equation}
\begin{lemme}
The system (\ref{3systeme3}) is global exact controllable implies that there exists a solution to the following controllability problem:\\
For any $h_{0},h_{1}\in H^{2}(0,1)$ with $0<c\leq h_{0}<M$, $0<c_{1}\leq h_{1}<M_{1}$, there exists control $v_{4}, v_{5}\in H^{1}(0,T)$ and a solution of (\ref{3systeme3}) such that $h(T,x)=h_{1}(x)$ in $(0,1)$ with $h\in Y_{T}=L^{2}(0,T;H^{3}(0,1))\cap H^{1}(0,T;H^{1}(0,1))$.
\end{lemme}
{\bf Proof:} By multiplying (\ref{3systeme2}) by $h$, we check that $h$ verifies the second equation of  (\ref{3systeme3}). In order to ensure the compatibility between the first and the second equation of (\ref{3systeme2}), we assume that:
\begin{equation}
v_{5}(t)=-\int^{t}_{0}v_{3}(s)ds.
\end{equation}
%ATTENTION ICI ON RESTREINT NOTRE CHOIX CE QUI N'EST PAS BON!!!! EN FAIT IL FAUT FAIRE LE DERNIER CHANGEMENT. $\blacksquare$\\
%\\
To summarize we have proved that if the system (\ref{3systeme}) is global exact controllable, it implies that the following heat equation:
\begin{equation}
\begin{cases}
\begin{aligned}
&\frac{\p}{\p t}h-\mu\p_{xx}h=0,\;\;(t,x)\in Q=(0,T)\times (0,1)\\
&h(t,0)=v_{6}(t),\;\;h(t,1)=v_{7}(t),\;\;t\in(0,T),\\
&h(0,x)=h_{0}(x).\\
\end{aligned}
\end{cases}
\label{3systeme4}
\end{equation}
is global exact controllable with:\\
$h_{0},h_{1}\in H^{2}(0,1)$ with $0<c\leq h_{0}<M$, $0<c_{1}\leq h_{1}<M_{1}$ and with $v_{6}=e^{v_{5}(t)},v_{7}=e^{v_{4}(t)}\in H^{1}(0,T)$ strictly positive where:
$$
\begin{aligned}
&h(t,x)=\exp{\big(-\frac{1}{\mu}\int^{x}_{0}u(t,y)dy-\int^{t}_{0}v_{3}(s)ds\big)},\\
&v_{3}(t)=-\frac{v_{1}^{2}(t)}{\mu}+u_{x}(t,0),
\end{aligned}
$$
and:
$$h_{1}(x)=K\exp{\big(-\frac{1}{\mu}\int^{x}_{0}u(t,y)dy}\big),$$
with $K=-\int^{T}_{0}v_{3}(s)ds$.
%The next step is to prove that the previous controllability result is equivalent to a controllability problem for a linear heat equation with two time-dependent controls, one of bilinear nature (multiplying the state function) and the other one acting at $x = 1$. By multiplying the previous system by $h$, we have:
%\begin{equation}
%\begin{cases}
%\begin{aligned}
%&\frac{\p}{\p t}z-\mu\p_{xx}z=0,\;\;(t,x)\in Q=(0,T)\times (0,1)\\
%&z(t,0)=e^{v_{5}(t)},\;\;z(t,1)=e^{v_{4}(t)},\;\;t\in(0,T),\\
%&z(0,x)=z_{0}(x).\\
%\end{aligned}
%\end{cases}
%\label{3systeme3}
%\end{equation}
%To do this we set:
%$$z(t,x)= \exp(-\int^{t}_{0}(v^{'}_{5}(s)+v_{3}(s))ds)h(t,x)= \exp(-v_{5}(t)-\int^{t}_{0}v_{3}(s))ds)h(t,x).$$
%We have seen that our system is equivalent to proving the exact controllability problem of the following system:
%\begin{equation}
%\begin{cases}
%\begin{aligned}
%&\frac{\p}{\p t}z-\mu\p_{xx}z=0,\;\;(t,x)\in Q=(0,T)\times (0,1)\\
%&z(t,0)=v_{7}(t),\;\;z(t,1)=v_{8}(t),\;\;t\in(0,T),\\
%&z(0,x)=z_{0}(x).\\
%\end{aligned}
%\end{cases}
%\label{3systeme4}
%\end{equation}
%More precisely for any $h_{0}, h_{1}\in H^{2}(0,1)$ with $h_{0}(x), h_{1}(x)>0$ in $(0,1)$ and $h_{0}(0)=h_{1}(0)=1$, there exists a constant $K>0$ and two controls $v_{7}(t)$, $v_{8}(t)\in H^{1}(0,T)$ which are strictly positive in $[0,T]$ such that the solution of (\ref{3systeme3}) satisfies $h(T,x)=K h_{1}(x)$ in $(0,1)$.
\section{ Proofs of theorem \ref{theo1} and \ref{theo2}}
\label{section3}
In this section, we will prove both theorems \ref{theo1} and \ref{theo2}.
We start with giving  technical results for the dual system of the heat equations that we shall use for proving theorem \ref{theo1} and \ref{theo2}:
\begin{proposition}
\label{lemma4}
 Let $0<\xi_{0} <\xi_{1} <\xi_{2} <1$. Then, for each $\theta>0$ there exists a time $T =T(\theta)>0 $such that the solution of the backwards heat equation:
 \begin{equation}
 \begin{cases}
&-U_{t}-U_{xx} =0,\;\; (t,x)\in(0,T^{*})\times(0,1),\\
&U(t,0)=U(t,1)=0,\;\; t\in (0,T^{*}), \\
&U(T^{*},x)=\delta_{\xi_{0}}-\theta  \delta_{\xi_{1}} +  \delta_{\xi_{2}},\;\;x\in(0,1)
\end{cases}
\label{21}
\end{equation}
satisfies
$$U_{x}(t,0)>0\;\;\mbox{and}\;\; U_{x}(t,1)<0\;\;\forall t\in(0,T^{*}).$$
Here $\delta_{x}$ is the Dirac measure at $x$.
\end{proposition}
%In the previous lemma, we have denoted by ?x the Dirac mass distribution at point x.
The following proposition is a generalization of the previous one  which was proved in \cite{2}:
\begin{proposition}
\label{lemma5}
 Let $T > 0$ and  $0<\xi_{0} <\xi_{1} <\xi_{2} <1$. Then there exists $\theta>0$ such that the solution of the backwards heat equation:
 \begin{equation}
 \begin{cases}
&-U_{t}-U_{xx} =0,\;\; (t,x)\in(0,T^{*})\times(0,1),\\
&U(t,0)=U(t,1)=0,\;\; t\in (0,T^{*}), \\
&U(T^{*},x)=\delta_{\xi_{0}}-\theta  \delta_{\xi_{1}} +  \delta_{\xi_{2}},\;\;x\in(0,1)
\end{cases}
\label{21}
\end{equation}
satisfies
$$U_{x}(t,0)>0\;\;\mbox{and}\;\; U_{x}(t,1)<0\;\;\forall t\in(0,T^{*}).$$
\end{proposition}
\subsection{Proof of theorem \ref{theo1} via the method of Guerrero et Imanuvilov (\cite{GI})}
We have verified that the global exact controllability of (\ref{3systeme}) implies the exact controllability problem of the heat equation:
\begin{equation}
\begin{cases}
\begin{aligned}
&\frac{\p}{\p t}h-\mu\p_{xx}h=0,\;\;(t,x)\in Q=(0,T)\times (0,1)\\
&z(t,0)=v_{7}(t),\;\;z(t,1)=v_{8}(t),\;\;t\in(0,T),\\
&h(0,x)=h_{0}(x).\\
\end{aligned}
\end{cases}
\label{3systeme4}
\end{equation}
 for any $h_{0}, h_{1}\in H^{2}(0,1)$ with $h_{0}(x), h_{1}(x)>0$ in $(0,1)$ and $h_{0}(0)=1$ and $h_{1}(0)=K>0$, and with two controls $v_{7}(t)$, $v_{8}(t)\in H^{1}(0,T)$ which are strictly positive in $[0,T]$ such that the solution of (\ref{3systeme3}) satisfies $h(T,x)= h_{1}(x)$ in $(0,1)$.\\
%More precisely, we have proved (see Lemmas 1Ð3) that the null controllability property for system (\ref{3systeme}) is equivalent
%to the existence of a positive constant K and positive controls v?1, v?2 ? H 1(0, T ) such that the solution h ? X1 of (23) satisfies h(T,x)=K in(0,1).\\
We are going to show the theorem \ref{theo1} by contradiction. Thus, we assume  that for any $T > 0$ and any $h_{0}\in H^{2}(0, 1)$ with $h_{0} >0$ and $h_{0}(0)=1$, there exists a constant $K >0$ and two controls $0<v_{1}\in H^{1}(0,T)$ and $0<v_{2}\in H^{1}(0,T)$ such that the solution of (\ref{3systeme4}) satisfies
$h(T,x)=K$ in $(0,1)$.\\
%Then, let us consider the function U given by proposition \ref{lemma4} for some $\theta\geq  2$, which is thus defined up to a time $T=T^{*}$. 
Multiplying the equation of (\ref{3systeme4}) by $U$ given by proposition \ref{lemma4} (for $\theta\geq 2$) and integrating over $(0, T)\times(0, 1)$ with $T=T^{*}$, we obtain that:
\begin{equation}
\int^{T}_{0}(U_{x}(t,0)\widetilde{v}_{1}(t)-U_{x}(t,1)\widetilde{v}_{2}(t)) dt +K(2-\theta) -\int^{1}_{0}U(0,x)h_{0}(x)dx =0
\label{24}
\end{equation}
for any $h_{0}\in H^{2} (0, 1)$. By using the facts that the normal derivative of $U$ is negative and $\theta\geq 2$, we observe that the two first terms of (\ref{24}) are non-positive.\\
Our goal now consists in choosing an initial condition $h_{0}$ with $h_{0}(0) = 1$ such that:
$$-\int^{1}_{0}U(0,x)h_{0}(x)dx<0,$$
in order to obtain a contradiction with the equality (\ref{24}) and the fact that the two first terms of this equality are non positive.\\
Using the fact that the normal derivative of $U$ is negative and $U$ verifies homogeneous Dirichlet boundary conditions, we observe by Taylor formula that there exists $\delta> 0$ such that:
$$U(0,x)\geq\delta x,\;\forall x\in(0,\delta)\;\;\mbox{and}\;\;U(0,x)\geq\delta(1-x),\;\forall x\in(1-\delta,1). $$
On the other hand since $U$ satisfies a heat equation (\ref{21}), we obtain an instantaneously regularizing effect on $U$ in Sobolev norms which depends on some power of $\frac{1}{T^{\alpha}}$ with $\alpha>1$,more precisely there exists a positive constant $C^{*}(T^{*},\theta)$ such that:
\begin{equation}
\|U(0,\cdot)\|_{L^{2}(0,1)}\leq C^{*}.
\label{U}
\end{equation}
We are going to follow the choice of $h_{0}$ as in \cite{GI}, let $h_{0} = h_{0}(x)\in (0,1)$ for all $x\in(0,1)$ such that 
$$h_{0}(x)=\frac{4 C^{*}}{\delta^{3}}\;\;\forall x\in(\frac{\delta}{4},\frac{3\delta}{4})U(1-\frac{3\delta}{4},1-\frac{\delta}{4}) ,$$
we have by using the fact that $h_{0} $ is positive:
$$
\begin{aligned}
-\int^{1}_{0}U(0,x)h_{0}(x)dx\leq&-\int^{\frac{3\delta}{4}}_{\frac{\delta}{4}}\frac{4 C^{*}}{\delta^{3}}U(0,x)dx-\int^{1-\frac{\delta}{4}}_{1-\frac{3\delta}{4}}\frac{4 C^{*}}{\delta^{3}}U(0,x)dx-\int^{1-\delta}_{0}U(0,x)h_{0}(x)dx\\
\leq&-C^{*}+(1-2\delta)C^{*}=-2\delta C{*}<0.
\end{aligned}
$$
This achieves the proof. $\blacksquare$
\subsection{Proof of theorem \ref{theo2}}
In this paragraph, we prove the theorem \ref{theo2}. The proof follows the same lines than in \cite{GI} and than in the previous section. That why we are just pointing out the differences and explaining how to deal with them.\\% Going back again to the previous section, we consider the following control system:
%We have seen that our system is equivalent to proving the exact controllability problem of the following system:
%\begin{equation}
%\begin{cases}
%\begin{aligned}
%&\frac{\p}{\p t}z-\mu\p_{xx}z=0,\;\;(t,x)\in Q=(0,T)\times (0,1)\\
%&z(t,0)=v_{7}(t),\;\;z(t,1)=v_{8}(t),\;\;t\in(0,T),\\
%&z(0,x)=z_{0}(x).\\
%\end{aligned}
%\end{cases}
%\label{3systeme4}
%\end{equation}
We are interested in showing that for any $h_{0}, h_{1}\in H^{2}(0,1)$ with $h_{0}(x), h_{1}(x)>0$ in $(0,1)$ and $h_{0}(0)=1$, $h_{1}(0)=K$, there exists two controls $v_{7}(t)$, $v_{8}(t)\in H^{1}(0,T)$ which are strictly positive in $[0,T]$ such that the solution of (\ref{3systeme3}) satisfies $h(T,x)= h_{1}(x)$ in $(0,1)$.\\
As in the previous section, we multiply the system (\ref{3systeme4}) by $U$ checking the proposition \ref{lemma5}, we have then:
\begin{equation}
\int^{T^{*}}_{0}(U_{x}(t,0)\widetilde{v}_{1}(t)-U_{x}(t,1)\widetilde{v}_{2}(t)) dt +\big(h_{1}(\xi_{0})-\theta h_{1}(\xi_{1})+h_{1}(\xi_{2})\big) =0
\label{24a}
\end{equation}
for all $0<h_{0}\in H^{2} (0, 1)$ and $0<h_{1}\in H^{2} (0, 1)$. \\
The first term of (\ref{24a}) is non positive, we are then using the same argument of contradiction than in the previous proof and in particular showing that:
$$-\int^{1}_{0}U(0,x)h_{0}(x)dx<0,$$
for a suitable choice on $h_{0}$. We proceed exactly as in the previous section with the same $h_{0}$ and furtheremore we take $h_{1}$ such that:
$$h_{1}(\xi_{0})-\theta h_{1}(\xi_{1})+h_{1}(\xi_{2})<0.$$ 
This is a contradiction with identity (\ref{24a}) which concludes the proof. $\blacksquare$%\\
 \\
 \\
\textbf{Aknowledgements}: The authors would like to thank Jean-Michel Coron for many helpful and fruitful discussions.

\end{document}